\documentclass[11pt]{article}
\title{\LARGE  Series  expansion for $L^p$  Hardy inequalities} 

\author{\Large G. Barbatis$^{ 1}$ ~~~~
 S. Filippas$^{ 1}$ ~~\& ~~ A. Tertikas$^{2,3}$  \\
                                                                           \\
        Department of Applied Mathematics$^{1}$ \\
         University of Crete,
         71409 Heraklion,  Greece \\ 
      gbarbati@tem.uoc.gr,~~  filippas@tem.uoc.gr\\ 
                                          \\    
 Department of Mathematics$^{2}$ \\
         University of Crete,
         71409 Heraklion,  Greece \\ 
          tertikas@tem.uoc.gr\\
                                      \\
        Institute of Applied and Computational Mathematics$^3$, \\
        FORTH, 71110 Heraklion, Greece \\ 
    \\ }

\begin{document}

\newcommand{\ana}{\nabla}
\newcommand{\R}{I \!  \! R}

\maketitle

\newcommand{\be}{\begin{equation}}
\newcommand{\ee}{\end{equation}}
\newcommand{\bea}{\begin{eqnarray}}
\newcommand{\eea}{\end{eqnarray}}
\newcommand{\la}{\label}
\newcommand{\xa}{\alpha}
\newcommand{\xb}{\beta}
\newcommand{\xg}{\gamma}
\newcommand{\xG}{\Gamma}
\newcommand{\xd}{\delta}
\newcommand{\xD}{\Delta}
\newcommand{\xe}{\varepsilon}
\newcommand{\xz}{\zeta}
\newcommand{\xh}{\eta}
\newcommand{\Th}{\Theta}
\newcommand{\xk}{\kappa}
\newcommand{\xl}{\lambda}
\newcommand{\xL}{\Lambda}
\newcommand{\xm}{\mu}
\newcommand{\xn}{\nu}
\newcommand{\ks}{\xi}
\newcommand{\KS}{\Xi}
\newcommand{\xp}{\pi}
\newcommand{\xP}{\Pi}
\newcommand{\xr}{\rho}
\newcommand{\xs}{\sigma}
\newcommand{\xS}{\Sigma}
\newcommand{\xf}{\phi}
\newcommand{\xF}{\Phi}
\newcommand{\ps}{\psi}
\newcommand{\PS}{\Psi}
\newcommand{\xo}{\omega}
\newcommand{\xO}{\Omega}
\newcommand{\Ren}{ I \! \! R^N}
\newcommand{\Real}{ I \! \! R}
\newcommand{\Br}{B_r}
\newcommand{\bBr}{\partial \! B_r}
\newcommand{\ra}{\rightarrow}
\newcommand{\rft}{\rightarrow +\infty}
\newcommand{\bin}{\int_{\bBr}} 
\newcounter{newsection}
\newtheorem{theorem}{Theorem}[section]
\newtheorem{lemma}[theorem]{Lemma}
\newtheorem{prop}[theorem]{Proposition}
\newtheorem{coro}[theorem]{Corollary}
\newtheorem{defin}[theorem]{Definition}
\newcounter{newsec} \renewcommand{\theequation}{\thesection.\arabic{equation}}

 \newcommand{\nl}{\newline}
 \newcommand{\dist}{{\rm dist}}
 \newcommand{\N}{{\bf N}}
 \newcommand{\Z}{{\bf Z}}
 \newcommand{\C}{{\bf (C)}}
 \newcommand{\hil}{{\cal H}}
 \newcommand{\cC}{{\cal C}}
 \newcommand{\cF}{{\cal F}}
 \newcommand{\cB}{{\cal B}}
 \newcommand{\cL}{{\cal L}}
 \newcommand{\cD}{{\cal D}}
 \newcommand{\cS}{{\cal S}}
 \newcommand{\cE}{{\cal E}}
 \newcommand{\cA}{{\cal A}}
 \newcommand{\CC}{{\bf (C')}}

 \newcommand{\re}{{\rm Re}\;}
 \newcommand{\im}{{\rm Im}\;}

 \newcommand{\diam}{{\rm diam}}
 \newcommand{\all}{\mbox{ all }}
 \newcommand{\as}{\mbox{ as }}
 \newcommand{\diver}{{\rm div}}
 \newcommand{\supp}{{\rm supp}}
 \newcommand{\inprod}[2]{{\langle{#1},{#2}\rangle}}
 \newcommand{\vol}{{\rm vol}}
 \newcommand{\codim}{{\rm codim}}
 \newcommand{\half}{\frac{1}{2}}
 \newcommand{\bb}{\hfill}
 \newcommand{\darr}[4]{{\left\{\begin{array}{ll}
   {#1}&{#2}\\
   {#3}&{#4}
 \end{array}\right.}}
 \newcommand{\tarr}[6]{{\left\{\begin{array}{lll}
   {#1}&{#2}\\
   {#3}&{#4}\\
   {#5}&{#6}
\end{array}\right.}}
\newcommand{\darrn}[4]{{\begin{array}{ll}
   {#1}&{#2}\\
   {#3}&{#4}
\end{array}}}
\newcommand{\tarrn}[6]{{\begin{array}{lll}
   {#1}&{#2}\\
   {#3}&{#4}\\
   {#5}&{#6}
\end{array}}}
\newcommand{\ia}{({\rm i})}
\newcommand{\ib}{({\rm ii})}
\newcommand{\ic}{({\rm iii})}
\newcommand{\hkp}{{\biggl|\frac{k-p}{p}\biggl|}}


 \parindent=0pt
 \parskip=5pt


\begin{abstract}

We consider a general class of sharp $L^p$ Hardy inequalities in $\R^N$ involving distance
from a surface of general codimension $1\leq k\leq N$.
We show that we can succesively improve them by adding to the right hand side
a lower order term with optimal weight and best constant. This leads to an infinite series
improvement of $L^p$ Hardy inequalities.

\end{abstract}

\section{Introduction}

Let $\xO$ be a bounded domain in $\Ren$ containing the origin. 
 Hardy  inequality  asserts that for  any $p>1$
 \begin{equation}
 \int_{\Omega}|\nabla u|^pdx\geq
 \bigg|\frac{N-p}{p}\bigg|^p\int_{\Omega }\frac{|u|^p}{|x|^p}dx,
 ~~~~~~~u\in C^{\infty}_{c}(\Omega \setminus\{0\}),
 \label{eq:h1}
 \end{equation}
 with  $|\frac{N-p}{p}|^p$ being the best constant,
see for example \cite{HLP}, \cite{OK}, \cite{DH}.
An analogous result  asserts that for a convex domain  $\Omega\subset \R^N$ with
 smooth boundary, and $d(x) ={\rm dist}(x,\partial\Omega)$, there holds
 \begin{equation}
 \int_{\Omega}|\nabla u|^pdx\geq
 \Bigl( \frac{p-1}{p} \Bigr)^p\int_{\Omega}\frac{|u|^p}{d^p}dx,
 ~~~~~~~u\in C^{\infty}_{c}(\Omega),
 \label{eq:h2}
 \end{equation}
 with $(\frac{p-1}{p})^p$ being the best constant, cf
 \cite{MS}, \cite{MMP}. See \cite{OK} for a comprehensive account of Hardy inequalities
and \cite{D} for a review of recent results.

In the last few years  improved versions of the above  inequalities have been obtained,
in the sense that nonnegative  terms  are added in the right hand side of (\ref{eq:h1})
or (\ref{eq:h2}). Improved Hardy inequalities are useful in the study of critical
phenomena in elliptic and parabolic PDE's; see, e.g., [BM, BV, MMP, VZ]. 
 In this work we obtain an infinite series improvement for general
Hardy inequalities, that include  (\ref{eq:h1}) or (\ref{eq:h2}) as special cases.

Before stating our main theorems  let us first introduce some notation.
Let $\xO$   be a domain in $\R^N$,   $N \geq 2$,  and $K$ 
a piecewise smooth surface of codimension $k$, $k=1,\ldots,N$. In case $k=N$,
we adopt the convention that $K$ is a point, say, the origin.
 We also  set
\[d(x) = {\rm dist}(x, K),
\]
  and we assume that the
 following inequality
holds  in the weak sense:
\[p \neq k, \qquad \qquad \Delta_p d^{\frac{p-k}{p-1}}\leq 0, 
\quad  \quad \mbox{ in $\Omega  \setminus K$.} \qquad \hfill \hspace{2cm}  \C
\]
Here $\Delta_p$ denotes the usual $p$-Laplace operator,
$\Delta_p w=\diver (|\nabla w|^{p-2}\nabla w)$. 
When $k=N$  (C) is  satisfied as equality, since $d^{\frac{p-k}{p-1}}=
|x|^{\frac{p-N}{p-1}}$ is the fundamental solution of the $p$-Laplacian.
Also, if $k=1$, $\Omega$ is convex and  $K=\partial\Omega$ condition (C)
is satisfied.  For a detailed analysis of this condition, as well as
for  examples in the intermediate cases  $1<k<N$,  we refer to [BFT].

We next  define the function:
\begin{equation}
X_1(t) = (1- \log t)^{-1},\quad t\in (0,1),
\label{eq:flash}
\end{equation}
and recursively
\[
X_k(t) = X_{1}(X_{k-1}(t)), ~~~~~~~k=2,3,\ldots ;
\]
these are iterrated logarithmic functions suitably normalized.
We also set
\be
I_m[u] := \int_{\xO} |\nabla u|^p dx - |H|^p
 \int_{\xO} \frac{|u|^p}{d^p} dx -  \frac{p-1}{2p}|H|^{p-2}   \sum_{i=1}^{m}
\int_{\xO} \frac{|u|^p}{d^p}
 X_1^2 X_2^2
 \ldots X_i^2 dx.
\la{1.3}
\ee
where
\[
H = \frac{k-p}{p}.
\]
Our main result reads

{\bf Theorem A }
{\em
 Let $\Omega$ be a domain in $\R^N$ and $K$ 
a piecewise smooth surface of codimension $k$, $k=1,\ldots,N$.
Suppose that  $\sup_{x\in\Omega}d(x) < \infty$
and       condition (C) is satisfied. Then: \nl
(1) There exists a positive constant
$D_0=D_0(k,p)\geq \sup_{x\in\Omega}d(x)$ such that for any $D \geq  D_0$
and all  $u \in W^{1,p}_0(\Omega\setminus K)$ there holds
\bea
\int_{\Omega}|\nabla u|^pdx - |H|^p\int_{\Omega}\frac{|u|^p}{d^p}dx
\geq ~~~~~~~~~~~~~~~~~~~~~~~~~~~~~~~~~~~~~~~~~~~~
\nonumber \\
~~~~~~~~~~~~~ 
 \frac{p-1}{2p}|H|^{p-2} \left( \sum_{i=1}^{\infty}
 \int_{\Omega}\frac{|u|^p}{d^p}X_1^2(d/D) \ldots X_i^2(d/D)    dx \right).
\la{1.1}
\eea
If in addition $2 \leq p <k$, then we can take
$D_0 = \sup_{x \in \Omega} d(x)$. \nl
(2) Moreover, for each $m=1,2,\ldots$ the constant
 $\frac{p-1}{2p}|H|^{p-2}$ is the best constant for the 
corresponding $m$-Improved Hardy inequality, that is,
\[
  \frac{p-1}{2p}|H|^{p-2} = \inf_{u \in W^{1,p}_0(\Omega\setminus K) }
 \frac{I_{m-1}[u]  }
{\int_{\xO}  \frac{ |u|^p}{d^p} X_1^2 X_2^2
 \ldots X_m^2   dx},
\]
in  either of
the following cases: (a) $k=N$ and $K=\{0\} \subset \xO$, (b) $k=1$ and
 $K = \partial\Omega$, (c) $2 \leq k \leq N-1$ and $\xO \cap K \neq \emptyset$.
}

We also note that the exponent two of the logarithmic corrections
in (\ref{1.1}) are optimal; see Proposition \ref{thm:best} for the precise statement.

For $p=2$, $\xO$ convex and $K = \partial\Omega$  the first term in the infinite
series of (\ref{1.1}) was obtained in [BM]. In the more general framework
of Theorem A, the first term in the above series  was obtained in [BFT].
On the other hand, when $p=2$ and $K=\{0\}$ the full series was obtained
in [FT] by a different method. For other types of improved Hardy inequalities
we refer to \cite{BV, GGM, M, VZ}; in all these works one correction term is added
in the right hand side of the plain Hardy inequality.

We next consider the degenerate case $p=k$ for which we do not have 
the usual  Hardy inequality.
 In [BFT] a substitute for Hardy inequality was
given in that case.  The analogue of condition (C) is now:
\[p = k, \qquad \qquad   \Delta_p(-\ln  d) \leq 0, 
\quad  \quad \mbox{ in $\Omega  \setminus K$.} \qquad \hfill \hspace{2cm}  \CC
\]
If $(C')$ is satisfied then for any
$D\geq\sup_{\Omega}d(x)$ there holds (cf [BFT], Theorems 4.2 and 5.4):
\be
\int_{\Omega}|\nabla u|^k dx\geq \left(\frac{k-1}{k}\right)^k
\int_{\Omega}\frac{|u|^k}{d^k}X_1^k(d/D)dx,\quad u\in W^{1,p}_0
(\Omega\setminus K),
\label{1.5}
\ee
with $\left(\frac{k-1}{k}\right)^k$ being the best constant.
In our next result we obtain a series improvement  for  inequality (\ref{1.5}). We set
\bea
\tilde{I}_m[u] =: \int_{\Omega}|\nabla u|^k dx  -   \left(\frac{k-1}{k}\right)^k
\int_{\Omega}\frac{|u|^k}{d^k}X_1^k(d/D)dx   -
~~~~~~~~~~~~~~~~~~~~~~~~~~~ \nonumber \\
-\frac{1}{2} \left(\frac{k-1}{k}\right)^{k-1} \sum_{i=2}^{m} 
\int_{\Omega}\frac{|u|^k}{d^k}X_1^k(d/D)  X_2^2(d/D)
 \ldots X_m^2(d/D)   dx.  \nonumber 
\eea
We then have

{\bf Theorem B}
{\em
  Let $\Omega$ be a domain in $\R^N$ and $K$ 
a piecewise smooth surface of codimension $k$, $k=2,\ldots,N$.
 Suppose that  $\sup_{x\in\Omega}d(x) < \infty$
and   condition (C') is satisfied. Then,

(1) for any
$D\geq\sup_{\Omega}d(x)$  and all  $u \in W^{1,k}_0(\Omega\setminus K)$   there holds
\bea
\int_{\Omega}|\nabla u|^k dx  -   \left(\frac{k-1}{k}\right)^k
\int_{\Omega}\frac{|u|^k}{d^k}X_1^k(d/D)dx   \geq  
~~~~~~~~~~~~~~~~~~~~~~~~ \nonumber \\
\frac{1}{2} \left(\frac{k-1}{k}\right)^{k-1} \sum_{i=2}^{\infty} 
\int_{\Omega}\frac{|u|^k}{d^k}X_1^k(d/D)  X_2^2(d/D)
 \ldots X_i^2(d/D)   dx.
\la{1.6}
\eea
(2) Moreover, for each $m=2,3, \ldots$ the constant $\frac{1}{2} \left(\frac{k-1}{k}\right)^{k-1}$
is the best constant for the corresponding $m$-Improved  inequality. That is
\[
  \frac{1}{2} \left(\frac{k-1}{k}\right)^{k-1}    = \inf_{u \in W^{1,p}_0(\Omega\setminus K) }
 \frac{\tilde{I}_{m-1}[u]  }
{\int_{\xO}  \frac{ |u|^k}{d^k} X_1^k X_2^2
 \ldots X_m^2   dx}.
\]
in  either of
the following cases: (a) $k=N$ and $K=\{0\} \subset \xO$, (b) $2 \leq k \leq N-1$ and
 $\xO \cap K \neq \emptyset$.
}

To prove parts (1) of  the above Theorems, we make use of suitable vector fields and
elementary inequalities; this is carried out in Section 2. To prove the second parts,
we use a local argument and appropriate  test functions; this is done in Section 3.


\setcounter{equation}{0}
\section{The series expansion}
In  this Section we will derive the series improvement that appear in part (1)
of Theorems A and B.
We shall repeatedly use the differentiation rule
\begin{equation}
 \frac{d}{dt} X_i^{\beta}(t)=\frac{\beta}{t}X_1(t)X_2(t)\ldots X_{i-1}(t)X_i^{1+\beta}(t),\qquad i=1,2,\ldots,\quad \beta\neq -1,
\label{eq:3.52}
\end{equation}
which is proved by induction: for $i=1$ (\ref{eq:3.52}) follows immediately from the definition of $X_1(t)$, cf. (\ref{eq:flash}):
\[ \frac{d}{dt}X_1^{\beta}(t)=\frac{\beta}{t}(1-\log t)^{-\beta-1}=\frac{\beta}{r}X_1^{\beta+1}(t).\]
Moreover assuming (\ref{eq:3.52}) for a fixed $i\geq 1$ we have
\begin{eqnarray*}
\frac{d}{dt}X_{i+1}^{\beta}(t)&=&\frac{d}{dt}[X_1^{\beta}(X_i(t))]\\
&=&\frac{\beta}{X_i(t)}X_1^{\beta+1}(X_i(t))\frac{dX_i(t)}{dt}\\
&=&\frac{\beta}{X_i(t)}X_{i+1}^{\beta+1}(t)\frac{1}{t}X_1(t)\ldots X_{i-1}(t)X_i^2(t)\\
&=&\frac{\beta}{t}X_1(t)\ldots X_i(t)X_{i+1}^{\beta+1}(t);
\end{eqnarray*}
hence (\ref{eq:3.52}) is proved.

{\em Proof of Theorem A(1):}
We will make use of a suitable vectot field as in \cite{BFT}.
If $T$ is  a $C^1$ vector field in $\Omega$, then,
for any  $u \in C_c^{\infty}(\Omega \setminus K)$ we first 
integrate by parts and then use  H\"{o}lder's inequality to obtain
 \begin{eqnarray*}
 \int_{\Omega} {\rm div}\, T ~  |u|^p dx & = &
 -p \int_{\Omega} ( T \cdot \nabla u ) |u|^{p-2} u dx \\
 & \leq & p \left( \int_{\Omega} | \nabla u |^p  dx \right)^{\frac{1}{p}}
 \left( \int_{\Omega} |T|^{\frac{p}{p-1}} |u|^{p} dx  \right)^{\frac{p-1}{p}} \\
 & \leq &
  \int_{\Omega}  | \nabla u |^p dx  + (p-1)
 \int_{\Omega} |T|^{\frac{p}{p-1}} |u|^{p} dx.
 \end{eqnarray*}
We therefore arrive at
\be
\int_{\Omega}  | \nabla u |^p  dx  \geq
  \int_{\Omega} ( {\rm div}  T - (p-1)  |T|^{\frac{p}{p-1}} ) |u|^{p} dx.
\la{2.21}
\ee
For $m \geq 1$ we introduce the notation
\begin{eqnarray*}
 && \eta(t) =  \sum_{i=1}^{m} X_1(t) \ldots X_i(t),  \\
 && B(t)= \sum_{i=1}^{m} X_1^2(t) \ldots X_i^2(t),\qquad\quad t\in (0,1).
 \end{eqnarray*}
In view of (\ref{2.21}) in order to prove (\ref{1.1}) it is enough to
establish the following pointwise estimate: 
\be
 {\rm div}\, T   - (p-1)  |T|^{\frac{p}{p-1}} \geq
\frac{|H|^{p}}{d^p}\left( 1  + \frac{p-1}{2p H^2} B(d(x)/D) \right).
\la{2.22}
\ee
To proceed we now make a specific choice of $T$. We take
 \[
 T(x) = H  |H|^{p-2}   \frac{\nabla d(x)}{d^{p-1}(x)}
 \left( 1 +  \frac{p-1}{pH}  \eta(d(x)/D)
 + a \eta^2(d(x)/D) \right).
 \]
where $a$ is a free parameter to be chosen later. In any case $a$  will be
such that the quantity $1 +  \frac{p-1}{pH} \eta (d/D) +a \eta^2(d/D)$ is positive
on $\Omega$. Note that $T(x)$
is singular at $x \in K$, but since  $u \in C_c^{\infty}(\Omega \setminus K)$
all previous calculations are legitimate.

When computing ${\rm div} T$ we need to diffferentiate $\eta(d(x)/D)$. Recalling (\ref{eq:3.52}) a straightforward
calculation gives
\[
\eta^{'}(t) = \frac{1}{t} \left( X_1^2  + ( X_1^2 X_2 + X_1^2 X_2^2 )
 + \cdots +( X_1^2 X_2  \ldots X_m + \cdots + X_1^2  \ldots X_m^2)
\right),
\]
from which follows that
\be
t \eta^{'} (t) = \frac{1}{2}  B(t)
+ \frac{1}{2}  \eta^{2} (t).
\la{2.23}
\ee
On the other hand, observing that
\[
\Delta_p d^{\frac{p-k}{p-1}}=\frac{p-k}{p-1}\biggl|\frac{p-k}{p-1}
\biggr|^{p-2}d^{-k}(d\Delta d+(1-k)|\nabla d|^2),\\
\]
condition (C) implies
\be
(p-k)(d\Delta d+1-k)\leq 0.
\la{2.24}
\ee
For the sake of simplicity we henceforth omit the argument $d(x)/D$ from $\eta(d(x)/D)$ and
$B(d(x)/D)$.
Using (\ref{2.23}) and (\ref{2.24}) a straightforward calculation shows that
\be
 {\rm div}\, T  \geq  \frac{|H|^p}{d^p} \left( p + \frac{p-1}{H} \eta + pa \eta^2 
+ \frac{p-1}{2pH^2} (B+\eta^2) + \frac{a}{H} (B + \eta^2) \eta
  \right).
\la{2.25}
\ee
It then follows that (\ref{2.22}) will be established once we prove
the following inequality
\[
(p-1) +  \frac{p-1}{H} \eta + (pa +\frac{p-1}{2pH^2}) \eta^2
+ \frac{a}{H} B \eta  + \frac{a}{H} \eta^3 - (p-1) 
\left( 1 +  \frac{p-1}{pH}  \eta 
 + a \eta^2  \right)^{\frac{p}{p-1}} \geq 0,
\]
for all $x\in\Omega$. We set for convenience
\bea
f(B,\eta) & = & (p-1) +  \frac{p-1}{H} \eta + (pa +\frac{p-1}{2pH^2}) \eta^2
+ \frac{a}{H} B \eta  + \frac{a}{H} \eta^3, \nonumber \\
g(\eta) & = & \left( 1 +  \frac{p-1}{pH}  \eta 
 + a \eta^2  \right)^{\frac{p}{p-1}},  \nonumber 
\eea
and the required inequality is written as
\be
f(B,\eta) -(p-1) g(\eta) \geq 0.
\la{2.27}
\ee
When  $\eta=\eta(d(x)/D)>0$ is small, 
 the Taylor expansion of $g(\eta)$   about $\eta=0$, gives
\be
g(\eta) = 1 + \frac{1}{H} \eta + \frac{1}{2}  \left( \frac{2ap}{p-1} + 
\frac{1}{p H^2} \right) \eta^2 + \frac{1}{6} \left( \frac{6a}{(p-1)H}
+ \frac{2-p}{p^2 H^3} \right) \eta^3 + O(\eta^4).
\la{2.28}
\ee
Let us also  note, that in the special case $a=0$,
there holds
\be
g(\eta) = 1  + \frac{1}{H} \eta  + 
\frac{1}{2p H^2} \eta^2  + \frac{2-p}{6p^2 H^3}  \left( 1 + \frac{p-1}{pH}
 \xi_{\eta} \right)^{\frac{3-2p}{p-1}}  \eta^3, ~~~~~~~(a=0),
\la{2.29}
\ee
for some $\xi_{\eta} \in (0, \eta)$, without any smallness assumption on
$\eta$.

In view of (\ref{2.28}), if $\eta$ is small,
inequality (\ref{2.27}) will be proved once we show:
\be
\frac{a}{H}   \geq     \left(   \frac{(2-p)(p-1)}{6 p^2 H^3} + O(\eta)
 \right) \frac{\eta^2}{B}.
\la{2.30}
\ee
From the definition of $\eta$ and $B$ it follows easily that
\be
m \geq \frac{\eta^2}{B} \geq 1.
\la{2.32}
\ee

We will show that  for any choice of $H$ and $p>1$, there exists an $a \in \R$,
such that (\ref{2.27}) holds true. 
We distinguish various cases:

\noindent {\bf (a)}  $H>0$, $1<p<2$.
We assume that $\eta$ is small, which amounts to taking $D$ big.
It is enough to show that we can choose $a$ such that  (\ref{2.30}) holds.
 In view of  (\ref{2.32}) we see that for
(\ref{2.30}) to be valid, it is enough to take $a$ to be big and positive.

\noindent {\bf (b)} $H>0$, $p \geq 2$. In this case we choose $a=0$ and we use (\ref{2.29}). Notice that
under our current assumptions on $H$, $p$ the last term in (\ref{2.29})
is negative and therefore
\be
g(\eta) \leq  1  + \frac{1}{H} \eta  + 
\frac{1}{2p H^2} \eta^2,~~~~~~~(a=0).
\la{2.34}
\ee
On the other hand
\[
f(B,\eta) =  (p-1) +  \frac{p-1}{H} \eta 
 +\frac{p-1}{2pH^2} \eta^2,~~~~~~~(a=0),
\]
and therefore (\ref{2.27}) is satisfied, without any smallness assumption
on $\eta$. In particular, we can take $D_0 = \sup_{x \in \Omega} d(x)$
  in this case.

\noindent {\bf (c)} $H<0$, $1<p<2$.
We assume that $\eta$ is small. In this case, the right hand side of
(\ref{2.30}) is negative. Hence, we can choose
$a=0$ and (\ref{2.30}) holds true.

\noindent {\bf (d)} $H<0$, $p \geq 2$. Arguing as in  case (a)
we take $a$ to be big and negative, and (\ref{2.30}) holds true. $\hfill //$

We next consider the degenerate case $p=k$.

{\em Proof of Theorem B(1):}
We  assume that $p=k \geq 2$ and that condition (C') is satisfied. The proof
is quite similar to the previous one.

An easy calculation shows that condition (C') implies that
\be
d \Delta d +1 -k \geq 0.
\la{2.35}
\ee

 We  now  choose the vector field
\be
T(x) = \Bigl( \frac{k-1}{k} \Bigr)^{k-1} \frac{\nabla d}{d^{k-1}} \Bigl(
X_1^{k-1} +
\sum_{i=2}^{m} X_1^{k-1}X_2 \ldots X_i \Bigr).
\ee
where, here and below, $X_j=X_j(d(x)/D)$. Taking into account (\ref{2.35}) a straightforward calculation yields that
\bea
{\rm div} T   - (p-1)  |T|^{\frac{p}{p-1}}   \geq  
 \frac{(k-1)^k}{k^{k-1}} \frac{ X_1^{k} }{d^k} 
 \Biggl(   1 + \sum_{i=2}^{m} X_2 \ldots X_i 
~~~~~~~~~~~~~~~~~~~~~~~~
 \nonumber \\
  ~~~~ + \frac{1}{k-1}
  \sum_{i=2}^m \sum_{j=2}^i X_2^2 \ldots X_j^2 X_{j+1} \ldots X_i 
  -  \frac{k-1}{k}  
\Bigl(
 1 + \sum_{i=2}^{m} X_2 \ldots X_i   \Bigr)^{\frac{k}{k-1}} \Biggr).
\label{2.39}
\eea
To estimate the last term in the right hand side of (\ref{2.39}) we use 
 Taylor's expansion to obtain the inequality
\[
\Bigl(
 1 + \sum_{i=2}^{m} X_2 \ldots X_i
 \Bigr)^{\frac{k}{k-1}}  
\leq 1 + \frac{k}{k-1} \sum_{i=2}^{m} X_2 \ldots X_i
+ \frac{k}{2(k-1)^2} \Bigl( \sum_{i=2}^{m} X_2 \ldots X_i \Bigr)^2.
\]
It then follows that
\bea
{\rm div} T   - (p-1)  |T|^{\frac{p}{p-1}}   \geq  
~~~~~~~~~~~~~~~~~~~~~~~~~~~~~~~~~~~~~~~~~~~~~~~~~~~~~~~~~~~~~~~~
\nonumber  \\
 \frac{(k-1)^{k-1}}{k^{k-1}} \frac{ X_1^{k} }{d^k} 
 \left( 
\frac{k-1}{k} +
  \sum_{i=2}^m \sum_{j=2}^i X_2^2 \ldots X_j^2 X_{j+1} \ldots X_i 
-   \frac{1}{2} 
\left( \sum_{i=2}^{m} X_2 \ldots X_i \right)^2  \right).
\la{2.40}
\eea
Expanding the square in the last term in (\ref{2.40}) we
conclude that
\[
{\rm div} T   - (p-1)  |T|^{\frac{p}{p-1}}   \geq  
\left( \frac{k-1}{k} \right)^{k-1} \frac{ X_1^{k} }{d^k} 
 \left( \frac{k-1}{k} + \frac{1}{2} \sum_{i=2}^m X_2^2 \ldots X_i^2 
\right),
\]
and the result follows. $\hfill //$

\setcounter{equation}{0}
\section{Best constants}

In this section we are going to prove the optimality of the Improved
Hardy Inequality of Section 2.
More precisely, for any $m\geq 1$ let us recall that
\begin{eqnarray*}
 I_m[u]&=&\int_{\Omega}|\nabla u|^pdx-|H|^p\int_{\Omega}\frac{|u|^p}{d^p}dx-\\
&&-\frac{p-1}{2p}|H|^{p-2}\int_{\Omega}\frac{|u|^p}{d^p}\left( X_1^2+X_1^2X_2^2+\cdots
+X_1^2\ldots X_m^2\right)dx,
\end{eqnarray*}
where $X_i=X_i(d(x)/D)$. We have the following

\begin{prop}
Let $\Omega$ be a domain in $\R^N$. $\ia$ If $2\leq k\leq N-1$ then we
take $K$ to be a piecewise smooth surface of codimension $k$ and assume $K\cap\Omega\neq 
\emptyset$; $\ib$ if $k=N$ then we take $K=\{0\}\subset\Omega$; $\ic$ if $k=1$ then we assume
$K =\partial \Omega$. Let $D\geq \sup_{\Omega} d(x)$ be fixed and
suppose that for some constants $B> 0$ and  $\gamma\in\R$ the following inequality
holds true for all $u\in W^{1,p}_0(\Omega\setminus K)$
\begin{equation}
I_{m-1}[u]   \geq  B \int_{\Omega}\frac{|u|^p}{d^p}X_1^2(d/D)\ldots 
X_{m-1}^2(d/D)X_m^{\gamma}(d/D)dx.
\label{eq:3.1}
\end{equation}
Then

$\ia$ $\gamma\geq 2$

$\ib$ If $\gamma = 2$ then $B  \leq \frac{p-1}{2p}|H|^{p-2}$.
\label{thm:best}
\end{prop}
{\em Proof.}
All our analysis will be local, say, in a fixed ball of radius  $\delta$
(denoted by $B_{\delta}$) centered at the origin, for some fixed small $\delta$.
The proof we present works for any $k=1,2,\ldots,N$. We note however that
for $k=N$ (distance from a point) the subsequent calculations are substantially
simplified, whereas for $k=1$ (distance from the boundary) one should
replace $B_{\delta}$ by $B_{\delta} \cap \Omega$. This last  change
entails some minor modifications, the arguments otherwise being  the same.
Without any loss of generality we may
assume that $0\in K    \cap \Omega$ ($k \neq 1$), or $0 \in \partial \Omega$ if $k=1$.
We divide the proof into several steps. 

{\bf Step 1. }Let $\phi\in C^{\infty}_c(B_{\delta})$ be such that
$0\leq \phi\leq 1$ in $B_{\delta}$ and
$\phi=1$ in $B_{\delta/2}$.
We fix small parameters $\alpha_0,\alpha_1,\ldots,\alpha_m>0$ and define the functions
\[ w(x)= d^{-H+\frac{\alpha_0}{p}}X_1^{\frac{-1+\alpha_1}{p}}(d/D)\ldots 
X_m^{\frac{-1+\alpha_m}{p}}(d/D)\]
and
\[ u(x)=\phi(x)w(x).\]
It is an immediate consequence of (\ref{eq:3.7}) below that $u\in W^{1,p}(\Omega)$ .
Moreover, if $k<p$ then $H<0$ and therefore $u|_K=0$.
On the other hand, if $k>p$ then a standard approximation
argument -- using cut-off functions -- shows that $W^{1,p}_0(\Omega\setminus K)=W^{1,p}(\Omega\setminus K)$.
Hence $u\in W^{1,p}_0(\Omega\setminus K)$. To prove the proposition we shall
estimate the corresponding Rayleigh quotient of $u$ in the limit
 $\alpha_0 \to 0$,      $\alpha_1 \to 0$, $\ldots$, $\alpha_m \to 0$ in this order.

It is easily seen that
\begin{equation}
\nabla w =-d^{\frac{-k+\alpha_0}{p}} X_1^{\frac{-1+\alpha_1}{p}}\!\!\!
\ldots X_m^{\frac{-1+\alpha_m}{p}}\Bigl(H+\frac{\zeta(x)}{p}\Bigr)\nabla d.
\label{eq:3.3}
\end{equation}
where
\begin{equation}
\zeta(x)=-\alpha_0 +(1-\alpha_1)X_1
+\ldots +(1-\alpha_m)X_1 X_2\ldots X_m,
\label{eq:3c}
\end{equation}
where, here and below, we omit the argument $d(x)/D$ from $X_i(d/D)$. Since $\xd$ is small the $X_i$'s are also small.  Hence
$\zeta(x)$ can be thought as a small parameter in the rest of the proof.

Now $\nabla u=\phi\,\nabla w+\nabla\phi\, w$ and
hence, using the elementary inequality
\begin{equation}
|a+b|^p  \leq |a|^p + c_p (|a|^{p-1} |b| + |b|^p), ~~~a,b \in \R^N, ~~~~~p>1,
\label{eq:3.2}
\end{equation}
we obtain
\begin{eqnarray}
\int_{\Omega}|\nabla u|^pdx 
 &  \leq  & \int_{B_{\delta}}\phi^{p} |\nabla w|^pdx
+c_p \,  \int_{B_{\delta}}  | \nabla \phi| |\phi|^{p-1}
 |\nabla  w|^{p-1}|w|\, dx  +\nonumber\\
&&+ c_p \int_{B_{\delta}} | \nabla \phi|^p |w|^p\, dx\nonumber  \\
& =: & I_1 + I_2 + I_3.  \label{eq:3.a}
\end{eqnarray}
We claim that
\begin{equation}
I_2,~I_3  = O(1) \quad \mbox{uniformly as $\alpha_0,\alpha_1,\ldots,\alpha_m$ tend to zero.}
\label{eq:3.4}
\end{equation}
Let us give the proof for $I_2$. Using the definition of $w(x)$ and
the regularity of $\phi$ we obtain
\[
I_2 \leq c\int_{B_{\delta}} d^{1-k +\alpha_0 } X_1^{-1+\alpha_1}\ldots
X_m^{-1+\alpha_m}\Bigl|H+\frac{\zeta(x)}{p}\Bigr|^{p-1}dx.
\] 
The appearance of $d^{-k+1}$ together with
the fact that $\zeta$ is small compared to $H$ implies 
that $I_2$ is uniformly bounded (see step 2). The integral $I_3$ is treated similarly.

{\bf Step 2.} We shall repeatedly  deal with integrals of the form
\begin{equation}
Q=\int_{\Omega}\phi^p d^{-k+\beta_0}X_1^{1+\beta_1}(d/D)\ldots X_m^{1+\beta_m}(d/D)dx,\quad \beta_i\in\R,
\label{eq:3.50}
\end{equation}
we therefore provide  precise conditions under which $Q<\infty$.
From our assumptions on $\phi$ we have
\[ \int_{B_{\delta/2}}d^{-k+\beta_0}X_1^{1+\beta_1}\!\!\ldots X_m^{1+\beta_m}dx
\leq Q\leq \int_{B_{\delta}}d^{-k+\beta_0}X_1^{1+\beta_1}\!\!\ldots X_m^{1+\beta_m}dx.\]
Using the coarea formula and the fact that
\[c_1r^{k-1}\leq \int_{\{d=r\}\cap B_{\delta}} dS <c_2r^{k-1}\]
we conclude that
\[c_1 \int_0^{\delta/2}r^{-1+\beta_0}X_1^{1+\beta_1}\!\!\ldots X_m^{1+\beta_m}dr
\leq Q\leq c_2\int_0^{\delta}r^{-1+\beta_0}X_1^{1+\beta_1}\!\!\ldots X_m^{1+\beta_m}dr.\]
where $X_i=X_i(r/D)$. Hence, recalling (\ref{eq:3.52}) we conclude that

\begin{equation}
Q<\infty \Longleftrightarrow \left\{
\begin{array}{ll}
& \beta_0>0 \\
\mbox{ or}&\mbox{$\beta_0=0$ and $\beta_1>0$}\\
\mbox{ or}&\mbox{$\beta_0=\beta_1=0$ and $\beta_2>0$}\\
& \cdots \\
\mbox{ or}&\mbox{$\beta_0=\beta_1=\ldots =\beta_{m-1}=0$ and $\beta_m>0$.}
\end{array}\right.
\label{eq:3.51}
\end{equation}

{\bf Step 3. }We introduce some auxiliary quantities and prove some
simple relations about them. For $0\leq i\leq j\leq m$ we define
\begin{eqnarray*}
&& A_0 =\int_{\Omega}\phi^p d^{-k+\alpha_0}X_1^{-1+\alpha_1}\!\!\!\!
\ldots X_m^{-1+\alpha_m} dx \\
&& A_i =\int_{\Omega}\phi^p d^{-k+\alpha_0}X_1^{1+\alpha_1}\!\!\!\!\ldots
X_i^{1+\alpha_i}X_{i+1}^{-1+\alpha_{i+1}}\!\!\!\!\ldots X_m^{-1+\alpha_m} dx \\
&& \Gamma_{0j}=\int_{\Omega}\phi^p d^{-k+\alpha_0}X_1^{\alpha_1}\!\!\!\!\ldots
X_i^{\alpha_i}X_{i+1}^{-1+\alpha_{i+1}}\!\!\!\!\ldots X_m^{-1+\alpha_m} dx \\
&& \Gamma_{ij}=\int_{\Omega}\phi^p d^{-k+\alpha_0}X_1^{1+\alpha_1}\!\!\!\!\ldots
X_i^{1+\alpha_i}X_{i+1}^{\alpha_{i+1}}\!\!\!\!\ldots X_j^{\alpha_j}
X_{j+1}^{-1+\alpha_{j+1}}\!\!\!\!\ldots X_m^{-1+\alpha_m} dx,
\end{eqnarray*}
with  $\Gamma_{ii}=A_i$. We have the following

{\em Two  identities: }
Let $0\leq i\leq m-1$ be given and assume that $\alpha_0 = \alpha_1 =\ldots=\alpha_{i-1}=0$. Then
\begin{eqnarray}
&&\alpha_i A_i=\sum_{j=i+1}^m (1-\alpha_j)\Gamma_{ij} + O(1)\label{eq:3.38}\\
&&\alpha_i\Gamma_{ij} =-\sum_{k=i+1}^j \alpha_k\Gamma_{kj} +\sum_{k=j+1}^m(1-\alpha_k)
\Gamma_{jk} +O(1)
\label{eq:3.39}
\end{eqnarray}
where the $O(1)$ is uniform as the $\alpha_i$'s tend to zero. Let us give the proof for
(\ref{eq:3.38}). We assume that $i>0$, the case $i=0$ being a straight-forward
adaptation. A direct computation gives
\begin{equation}
\alpha_i d^{-k}X_1\ldots X_{i-1}X_i^{1+\alpha_i}=\diver(d^{-k+1}X_i^{\alpha_i}\nabla d) -
d^{-k}(d\Delta d+1-k)X_i^{\alpha_i}.
\label{eq:3d}
\end{equation}

hence
\begin{eqnarray*}
\alpha_iA_i&=&\int_{\Omega}\phi^p\diver ( d^{-k+1}X_i^{\alpha_i}\nabla d)
X_{i+1}^{-1+\alpha_{i+1}}\ldots X_m^{-1+\alpha_m}dx-\\
&&-\int_{\Omega}\phi^p d^{-k} (d\Delta d+1-k)X_i^{\alpha_i}X_{i+1}^{-1+\alpha_{i+1}}\ldots
X_m^{-1+\alpha_m}dx \\
&=:&E_1-E_2.
\end{eqnarray*}
It is a direct consequence of [AS, Theorem 3.2] that
\begin{equation}
d\Delta d+1-k=O(d),\quad \mbox{as $d\to 0$,}
\label{eq:3e}
\end{equation}
hence $E_2$ is estimated by a constant
times $\int_{\Omega}\phi^p d^{-k+1}X_i^{\alpha_i}X_{i+1}^{-1+\alpha_{i+1}}\ldots X_m^{-1+\alpha_m}dx$
and therefore is bounded uniformly in $\alpha_0,\alpha_1,\ldots,\alpha_m$.
To handle $E_1$ we integrate by parts obtaining
\begin{eqnarray*}
E_1&=&-\int_{\Omega}\nabla\phi^p \cdot\nabla d\;d^{-k+1}X_i^{\alpha_i}X_{i+1}^{-1+\alpha_{i+1}}
\ldots X_m^{-1+\alpha_m}dx\\
&&-\int_{\Omega}\phi^p d^{-k+1}X_i^{\alpha_i}\nabla d\cdot\nabla\Bigl(X_1^{-1+\alpha_1}\ldots
X_m^{-1+\alpha_m}\Bigr)dx
\end{eqnarray*}
The first integral is of order $O(1)$ (similarly to $I_2,I_3$ above)
while the second is equal to $\sum_{j=i+1}^m(1-\alpha_j)\Gamma_{ij}$.
Hence (\ref{eq:3.38}) has been proved. To prove (\ref{eq:3.39}) we use (\ref{eq:3d}) once more and proceed
similarly; we omit the details.


{\bf Step 4. }We proceed to estimate $I_1$. It follows from (\ref{eq:3.3}) that
\[I_1=\int_{\Omega}\phi^p d^{-k+\alpha_0}X_1^{-1+\alpha_1}\ldots
X_m^{-1+\alpha_m}\Bigl|H+\frac{\zeta}{p}\Bigr|^pdx.\]
Since $\zeta$ is small compared to $H$ we may use Taylor's expansion  to obtain
\begin{equation}
\Bigl|H +\frac{\zeta}{p}\Bigr|^p \leq  |H|^p + |H|^{p-2} H \zeta + \frac{p-1}{2p} |H|^{p-2} \zeta^2
+ c |\zeta|^3.
\label{eq:3.16}
\end{equation}

Using this inequality we can bound $I_1$ by
\begin{equation}
I_1 \leq I_{10}+I_{11} + I_{12} + I_{13},
\label{eq:3.b}
\end{equation}
where
\begin{eqnarray}
I_{10}&=&|H|^p \int_{B_{\delta}} \phi^{p} 
 d^{-k +\alpha_0} X_1^{-1+\alpha_1}\!\!\!\ldots X_m^{-1+\alpha_m}dx=|H|^p\int_{\Omega}
\frac{|u|^p}{d^p}dx,\label{eq:30}\\ 
I_{11}&  = & |H|^{p-2} H 
 \int_{B_{\delta}} \phi^{p} 
 d^{-k +\alpha_0} X_1^{-1+\alpha_1}\!\!\!\ldots X_m^{-1+\alpha_m}
 \zeta(x) dx, \nonumber\\
I_{12}&  = & \frac{p-1}{2p}|H|^{p-2}
 \int_{B_{\delta}} \phi^{p} 
 d^{-k +\alpha_0} X_1^{-1+\alpha_1}\!\!\!\ldots X_m^{-1+\alpha_m}
  \zeta^2(x)dx, \nonumber\\
I_{13}&  = & c
 \int_{B_{\delta}} \phi^{p} 
 d^{-k +\alpha_0} X_1^{-1+\alpha_1}\!\!\!\ldots X_m^{-1+\alpha_m} 
|\zeta(x)|^3 dx.\nonumber
\end{eqnarray}

{\bf Step 5. }We shall prove that
\begin{equation}
I_{11},~I_{13}  = O(1) \quad \mbox{uniformly in $\alpha_0,\alpha_1,\ldots,\alpha_m$.}
\label{eq:3.31}
\end{equation}

Indeed, substituting for $\zeta$ in $I_{11}$ we see by a
direct application of (\ref{eq:3.38}) (for $i=0$) that $I_{11}=O(1)$.
To estimate $I_{13}$ we observe that $X_1\ldots X_i\leq cX_1$ for some $c>0$ and thus obtain
\begin{eqnarray*}
I_{13}&\leq& c_1\alpha_0^3 \int_{\Omega}\phi^pd^{-k+\alpha_0}X_1^{-1+\alpha_1}\!\!\!\ldots 
X_m^{-1+\alpha_m}dx+\\
&&+c_2\int_{\Omega}\phi^pd^{-k+\alpha_0}X_1^{2+\alpha_1}
X_2^{-1+\alpha_2}\!\!\!\!\ldots X_m^{-1+\alpha_m}dx.
\end{eqnarray*}

The second integral is bounded uniformly in the $\alpha_i$'s due to the factor $X_1^2$.
Moreover,
using the fact $0\leq\phi\leq 1$ and $\int_{\{d=r\}\cap B_{\delta}} dS <cr^{k-1}$ we obtain
\begin{eqnarray*}
&&\hspace{-4cm}\alpha_0^3\int_{\Omega}\phi^pd^{-k+\alpha_0}X_1^{-1+\alpha_1}\ldots X_m^{-1+\alpha_m}dx\\
&\leq&c\alpha_0^3\int_0^{\delta}r^{-1+\alpha_0}X_1^{-1+\alpha_1}(r/D)\ldots X_m^{-1+\alpha_m}(r/D)dr\\
&\leq&c\alpha_0^3\int_0^{\delta}r^{-1+\alpha_0}X_1^{-2}(r/D)dr \\
(r=Ds^{1/\alpha_0})\quad&=&cD^{\alpha_0}\alpha_0^2\int_0^{(\delta/D)^{\alpha_0}}
\Bigl(1-\frac{1}{\alpha_0}\log s\Bigr)^2 ds \\
&=& O(1)
\end{eqnarray*}
as $\alpha_0\to 0$, uniformly in $\alpha_1,\ldots\alpha_m$.
Hence (\ref{eq:3.31}) has been proved. Combining (\ref{eq:3.a}), (\ref{eq:3.4}),
(\ref{eq:3.b}), (\ref{eq:30}) and (\ref{eq:3.31}) we conclude that
\begin{equation}
\int_{\Omega}|\nabla u|^pdx -|H|^p\int_{\Omega}\frac{|u|^p}{d^p}dx
\leq I_{12}+O(1),
\label{eq:3.7}
\end{equation}
uniformly in the $\alpha_i$'s.

{\bf Step 6.} Recalling the definition of $I_{m-1}[\cdot]$ we obtain from (\ref{eq:3.7})
\begin{eqnarray}
I_{m-1}[u]&\leq&\frac{p-1}{2p}|H|^{p-2}
\int_{\Omega}\phi^pd^{-k+\alpha_0}X_1^{-1+\alpha_1}\!\!\!
\ldots X_m^{-1+\alpha_m}\times\nonumber\\
&&\hspace{3cm}\times\Bigl(\zeta^2(x)-\sum_{i=1}^{m-1}X_1^2\ldots X_i^2\Bigr)dx +O(1)\nonumber\\
&=:&\frac{p-1}{2p}|H|^{p-2}J +O(1).
\label{eq:3.8}
\end{eqnarray}
 Expanding  $\zeta^2(x)$ (cf (\ref{eq:3c})) and  collecting similar terms we obtain
\begin{eqnarray}
J&=&\int_{\Omega}\phi^pd^{-k+\alpha_0}X_1^{-1+\alpha_1}\ldots X_m^{-1+\alpha_m}
\left\{\alpha_0^2+\sum_{i=1}^m(1-\alpha_i)^2X_1^2\ldots X_i^2-\right.\nonumber\\
&&\hspace{2cm}-\sum_{i=1}^{m-1}X_1^2\ldots X_i^2 -2\alpha_0 \sum_{j=1}^m(1-\alpha_j)
X_1\ldots X_j +\nonumber\\
&&\left.\hspace{2cm}+2\sum_{i=1}^{m-1}\sum_{j=i+1}^{m}(1-\alpha_i)(1-\alpha_j)X_1^2\ldots X_i^2
X_{i+1}\ldots X_j\right\}dx\nonumber\\
&=&\alpha_0^2A_0+A_m +\sum_{i=1}^m(\alpha_i^2-2\alpha_i)A_i -2\alpha_0
\sum_{j=1}^m(1-\alpha_j)\Gamma_{0j}+\nonumber \\
&&\hspace{4cm}+\sum_{i=1}^{m-1}\sum_{j=i+1}^m
2(1-\alpha_i)(1-\alpha_j)\Gamma_{ij}.
\label{eq:3.40}
\end{eqnarray}

{\bf Step 7.} We intend to take the limit $\alpha_0\to 0$ in (\ref{eq:3.40}). All terms have
finite limits except those containing  $A_0$ and $\Gamma_{0j}$ which, when viewed separately,
 diverge. When combined
however, they give
\begin{eqnarray*}
&&\alpha_0^2A_0 -2\alpha_0
\sum_{j=1}^m(1-\alpha_j)\Gamma_{0j}\\
\mbox{(by (\ref{eq:3.38}))} &=&-\alpha_0 \sum_{j=1}^m(1-\alpha_j)\Gamma_{0j} +O(1)\\
\mbox{(by (\ref{eq:3.39}))} &=& -\sum_{j=1}^m(1-\alpha_j)\Bigl(
-\sum_{i=1}^j \alpha_i\Gamma_{ij} +\sum_{i=j+1}^m(1-\alpha_i)\Gamma_{ji} \Bigr)  +O(1)  \\
&=&\sum_{i=1}^m(\alpha_i-\alpha_i^2)A_i+\sum_{i=1}^{m-1}\sum_{j=i+1}^m(2\alpha_i-1)(1-\alpha_j)\Gamma_{ij} +O(1).
\end{eqnarray*}
All the terms in the last expression remain bounded as $\alpha_0\to 0$; hence taking the limit
in (\ref{eq:3.40}) we obtain
\begin{equation}
J=A_m-\sum_{i=1}^m\alpha_iA_i +\sum_{i=1}^{m-1}\sum_{j=i+1}^m(1-\alpha_j)\Gamma_{ij} +O(1)\qquad\quad
(\alpha_0=0)
\label{eq:3.41}
\end{equation}
where the $O(1)$ is uniform with respect to $\alpha_1,\ldots,\alpha_m$. 

{\bf Step 8.} We next let $\alpha_1\to 0$ in (\ref{eq:3.41}). All terms have finite
limits except those involving $A_1$ and $\Gamma_{1j}$ which diverge.
Using (\ref{eq:3.38}) once more -- this time for $i=1$ --
we see that when combined these terms stay bounded in the limit $\alpha_1\to 0$.
Hence
\begin{equation}
J=A_m-\sum_{i=2}^m\alpha_iA_i +\sum_{i=2}^{m-1}\sum_{j=i+1}^m(1-\alpha_j)\Gamma_{ij} +O(1)\qquad\quad
(\alpha_0=\alpha_1=0)
\label{eq:3.45}
\end{equation}

We proceed in this way and after letting $\alpha_{m-1}\to 0$ we are left with
\begin{equation}
J=(1-\alpha_m)A_m +O(1), ~~~~~~~~~~~~~~(\xa_0=\xa_1= \ldots \xa_{m-1}=0),
\label{eq:3.42}
\end{equation}
uniformly in $\alpha_m$.

Combining (\ref{eq:3.1}), (\ref{eq:3.8}) and (\ref{eq:3.42}) we conclude that
\begin{equation}
B\leq\frac{p-1}{2p}|H|^{p-2}\frac{(1-\alpha_m)A_m +O(1)}{\int_{\Omega}\phi^p 
d^{-k}X_1\ldots X_{m-1}X_m^{\gamma-1+\alpha_m}dx}.
\label{eq:3.43}
\end{equation}
Suppose now that $\gamma<2$. Then letting $\alpha_m\to 2-\gamma>0$ we observe
that the denominator in (\ref{eq:3.43}) tends to infinity while the numerator
stays bounded. This implies $B=0$ proving part (i) of the Proposition.

Now, if $\gamma=2$ then the denominator in (\ref{eq:3.43}) is equal to $A_m$.
Hence letting $\alpha_m\to 0$ we have $A_m\to\infty$ (by (\ref{eq:3.51})) and
hence $B\leq\frac{p-1}{2p}|H|^{p-2}$. This concludes the
proof. $\hfill //$



We next  consider the degenerate  case $p=k$.
We have the following
\begin{prop}
Let $\Omega$ be a domain in $\R^N$. $\ia$ If $2\leq k\leq N-1$ then we take $K$ to
be a piecewise smooth surface of codimension $k$ and assume $K\cap\Omega\neq 
\emptyset$; $\ib$ if $k=N$ then we take
$K=\{0\}\subset\Omega$.
Let $D\geq\sup_{x\in\Omega}d(x)$ be fixed and
suppose that for some constants $B> 0$ and  $\gamma\in\R$
the following inequality holds true for all $u\in C^{\infty}_{c}(\Omega\setminus K)$
\begin{equation}
\tilde{I}_{m-1}[u]   \geq  B \int_{\Omega}\frac{|u|^k}{d^k}X_1^k(d/D)X_2^2(d/D) \ldots X_m^{\gamma}(d/D)dx.
\label{eq:3.1a}
\end{equation}
Then:

$\ia$ $\gamma\geq 2$

$\ib$ If $\gamma = 2$ then $B  \leq \frac{1}{2}\bigl(\frac{k-1}{k}\bigr)^{k-1}$.
\label{thm:best1}
\end{prop}
{\em Proof. }The proof is similar to that of Proposition \ref{thm:best}.
Without any loss of generality we assume that $0\in K\cap \Omega$.
As in the previous theorem we let $\phi$ be a non-negative,
smooth cut-off function supported in $B_{\delta}=\{|x|<\delta\}$, equal
to one on $B_{\delta/2}$ and taking values in $[0,1]$.

Given small parameters $\alpha_1,\ldots,\alpha_m>0$ we define
\[w(x)=X_1^{\frac{-k+1+\alpha_1}{k}}(d/D)X_2^{\frac{-1+\alpha_2}{k}}(d/D)\ldots
X_m^{\frac{-1+\alpha_m}{k}}(d/D).\]
and
\[u(x)=\phi(x)w(x).\]
 Subsequent calculations will establish that $u\in W^{1,k}(\Omega)$ (see (\ref{eq:3.18})).
We will prove that $u\in W^{1,k}_0(\Omega\setminus K)$ by showing that
\begin{equation}
d^{\frac{\alpha_0}{k}}u\to u \quad\mbox{in $W^{1,k}(\Omega)$ as $\alpha_0\to 0$}.
\label{eq:3.60}
\end{equation}
We have
\begin{equation}
\int_{\Omega}|\nabla (d^{\frac{\alpha_0}{k}}u)-\nabla u|^k dx 
\leq c \alpha_0^{k}\int_{\Omega}d^{-k+\alpha_0}u^kdx
+\int_{\Omega}|d^{\frac{\alpha_0}{k}} -1 |^k |\nabla u|^kdx.
\label{eq:3.61}
\end{equation}
The second term in the right hand side of (\ref{eq:3.61}) tends to zero as $\alpha_0\to 0$
by the dominated convergence theorem. Moreover, there exists a constant $c_{\alpha_1}$
such that $X_1^{\alpha_1/2}X_2^{-1}\ldots X_m^{-1}\leq c_{\alpha_1}$. Hence the first
term in the right hand side of (\ref{eq:3.61}) is estimated by
\[ c_{\alpha_1}\alpha_0^k \int_{\Omega}d^{-k+\alpha_0}X_1^{-k+1+\frac{\alpha_1}{2}}dx.\]
A direct application of \cite[Lemma 5.2]{BFT} shows that this tends to zero
as $\alpha_0\to 0$. Hence $u\in W^{1,p}_0(\Omega\setminus K)$.

To proceed we use (\ref{eq:3.2}) obtaining
\begin{eqnarray}
\int_{\Omega}|\nabla u|^kdx 
 &  \leq  & \int_{\Omega} \phi^k |\nabla w|^kdx
+c_k \,  \int_{\Omega}  | \nabla \phi| |\phi|^{k-1}
 |\nabla  w|^{k-1}|w|\, dx  +\nonumber\\
&&+ c_k \int_{\Omega} | \nabla \phi|^k |w|^k\, dx\nonumber  \\
& =: & I_1 + I_2 + I_3.  \label{eq:3.15}
\end{eqnarray}
Arguing as in the proof of the previous proposition (cf step 1) we see that $I_2$ and $I_3$
are bounded uniformly with respect to the $\alpha_i$'s. Hence
\begin{equation}
\int_{\Omega}|\nabla u|^kdx \leq\int_{\Omega} \phi^k |\nabla w|^kdx +O(1)
\label{eq:3.17}
\end{equation}
uniformly as $\alpha_1,\ldots,\alpha_m\to 0$.

Now, a direct computation yields
\[
\nabla w=-d^{-1}X_1^{\frac{1+\alpha_1}{k}}X_2^{\frac{-1+\alpha_2}{k}}\!\!\!\! \ldots
X_m^{\frac{-1+\alpha_m}{k}}
\left(\frac{k-1}{k}+\frac{\zeta(x)}{k}\right) \nabla d
\]
where
\[\zeta(x)=-\alpha_1+\sum_{i=2}^m(1-\alpha_i)X_2(d/D)\ldots X_i(d/D).\]

From (\ref{eq:3.17}) and recalling (\ref{eq:3.16}) we have
\begin{eqnarray}
&&\hspace{-2cm}
\int_{\Omega}|\nabla u|^kdx\leq \int_{\Omega}\phi^k d^{-k}X_1^{1+\alpha_1}X_2^{-1+\alpha_2}
\!\!\!\!\ldots X_m^{-1+\alpha_m}\times \label{eq:3.18}\\
&&\times\left\{ \Bigl(\frac{k-1}{k}\Bigr)^k +\Bigl(\frac{k-1}{k}\Bigr)^{k-1}\zeta
+\frac{1}{2}\Bigl(\frac{k-1}{k}\Bigr)^{k-1}\zeta^2 +c|\zeta|^3\right\}dx.\nonumber
\end{eqnarray}
The term containing $|\zeta|^3$ is bounded uniformly with respect to
$\alpha_1,\ldots,\alpha_m$ (cf Step 4 in the previous proposition).
Moreover it is immediately seen that
\begin{equation}
\phi^k d^{-k}X_1^{1+\alpha_1}X_2^{-1+\alpha_2}
\!\!\!\!\ldots X_m^{-1+\alpha_m}=\frac{|u|^k}{d^k}X_1^k,
\label{eq:3aa}
\end{equation}
Hence
\begin{eqnarray}
\tilde{I}_{m-1}[u]
&\leq&\int_{\Omega}\phi^k d^{-k}X_1^{1+\alpha_1}X_2^{-1+\alpha_2}
\!\!\!\!\ldots X_m^{-1+\alpha_m}\times \nonumber\\
&&\qquad\left\{\Bigl(\frac{k-1}{k}\Bigr)^{k-1}(-\alpha_1+\sum_{i=2}^m(1-\alpha_i)X_2\ldots X_i)
+ \right. \nonumber\\
&&\qquad\left. +\frac{1}{2}\Bigl(\frac{k-1}{k}\Bigr)^{k-1}\Bigl(-\alpha_1+\sum_{i=2}^m(1-\alpha_i)
X_2\ldots X_i\Bigr)^2 -\right. \nonumber\\
&&\qquad\left.  - \frac{1}{2}\Bigl(\frac{k-1}{k}\Bigr)^{k-1}\sum_{i=2}^{m-1}
X_2^2\ldots X_i^2 \right\}dx +O(1)
\label{eq:3.19}
\end{eqnarray}
where the $O(1)$ is uniform with respect to all the $\alpha_i$'s. Expanding the square
and collecting similar terms we conclude that

\begin{equation}
\tilde{I}_{m-1}[u]\leq \frac{1}{2}\left(\frac{k-1}{k}\right)^{k-1}\tilde{J}+O(1),
\quad \mbox{uniformly in $\alpha_1,\ldots,\alpha_m$,}
\label{eq:3.35}
\end{equation}
where
\begin{equation}
\tilde{J}=A_m+\sum_{i=1}^m(\alpha_i^2-2\alpha_i)A_i+
2\sum_{i=1}^{m-1}\sum_{j=i+1}^m(1-\alpha_i)(1-\alpha_j)\Gamma_{ij}.
\label{eq:3.36}
\end{equation}

We intend to take the limit $\alpha_1\to 0$ in (\ref{eq:3.36}). All terms
have a finite limit except $A_1$ and $\Gamma_{1j}$ which do not contain the factor
$X_2^{1+\alpha_2}$. When combined they give

\begin{eqnarray*}
&&\hspace{-2cm}(\alpha_1^2-2\alpha_1)A_1 +2\sum_{j=2}^m(1-\alpha_1)(1-\alpha_j)\Gamma_{1j}\\
\mbox{(by (\ref{eq:3.38}))}&=&\alpha_1^2 A_1-2\alpha_1\sum_{j=2}^m(1-\alpha_j)
\Gamma_{1j} +O(1) \\
\mbox{(by (\ref{eq:3.38}))}&=&-\sum_{j=2}^m(1-\alpha_j)\alpha_1 \Gamma_{1j} +O(1) \\
\mbox{(by (\ref{eq:3.39}))}
&=&\sum_{j=2}^m(1-\alpha_j) \Bigl( \sum_{i=2}^j\alpha_i\Gamma_{ij} -\sum_{i=j+1}^m
(1-\alpha_i)\Gamma_{ji}\Bigr)
+O(1)\\
&=&\sum_{i=2}^m (\alpha_i-\alpha_i^2)A_i +
\sum_{i=2}^{m-1}\sum_{j=i+1}^m(2\alpha_i-1)(1-\alpha_j)\Gamma_{ij} +O(1).
\end{eqnarray*}
In this expression we can let $\xa_1 \ra 0$.
Hence (\ref{eq:3.36}) becomes
\begin{equation}
\tilde{J}=A_m -\sum_{i=2}^m\alpha_iA_i +\sum_{i=2}^{m-1}\sum_{j=i+1}^m (1-\alpha_j)\Gamma_{ij}
+O(1),
~~~~~~~~~~~~(\xa_1=0).
\label{eq:3.37}
\end{equation}
This relation is completely analogous to (\ref{eq:3.41}). For the rest of
the proof we argue as in the proof of Proposition \ref{thm:best};
we omit the details. $\hfill //$


\begin{thebibliography}{RRR}

\bibitem[AS]{AS}{Ambrosio L. and Soner H.M. Level set approach to mean curvature flow in arbitrary
codimension. J. Diff. Geometry {\bf 43} (1996) 693-737.}

\bibitem[BFT]{BFT}{ Barbatis G., Filippas S., and Tertikas A. A unified approach to improved
$L^p$ Hardy inequalities with best constants. submitted.}

 \bibitem[BM]{BM}{Brezis H. and Marcus M. Hardy's inequalities
 revisited. Ann. Scuola Norm. Pisa {\bf 25} (1997) 217-237.}

 \bibitem[BV]{BV}{Brezis H. and V\'{a}zquez J.-L. Blow-up solutions
 of some nonlinear elliptic problems. Rev. Mat. Univ. Comp. Madrid
 {\bf 10} (1997) 443-469.}

 \bibitem[D]{D}{Davies E.B. A review of Hardy inequalities. Oper. Theory Adv. Appl. {\bf 110}
(1998) 55-67.}

 \bibitem[DH]{DH}{Davies E.B. and Hinz A.M. Explicit constants for Rellich inequalities
 in $L_p(\Omega)$. Math. Z. {\bf 227} (1998) 511-523.}

 \bibitem[FT]{FT}{Filippas S. and Tertikas A. Optimizing Improved
 Hardy inequalities. J. Funct. Anal., to appear.}

\bibitem[GGM]{GGM}{Gazzola F., Grunau H.-Ch. and Mitidieri E. Hardy inequalities with optimal constants and
remainder terms. Trans. Amer. Math. Soc., to appear.}

 \bibitem[HLP]{HLP}{Hardy G., P\'{o}lya G. and Littlewood J.E. Inequalities. 2nd edition,
 Cambridge University Press 1952.}

 \bibitem[MMP]{MMP}{Marcus M., Mizel V.J. and Pinchover Y. On the best constant for Hardy's inequality
 in $\R^n$. Trans. Amer. Math. Soc. {\bf 350} (1998) 3237-3255.}

 \bibitem[MS]{MS}{Matskewich T. and Sobolevskii P.E. The best
 possible constant in generalized Hardy's inequality for
 convex domain in $\R^n$. Nonlinear Anal., Theory, Methods \& Appl.,
 {\bf 28} (1997) 1601-1610.}

 \bibitem[M]{M}{Maz'ja V.G. Sobolev spaces. Springer 1985.}

 \bibitem[OK]{OK}{Opic B. and Kufner A. Hardy-type inequalities. Pitman
 Research Notes in Math., vol.219, Longman 1990.}

 \bibitem[VZ]{VZ}{Vazquez J.L. and Zuazua E., The Hardy inequality and
 the asymptotic behavior of the heat equation with an
 inverse-square potential. J. Funct. Anal. {\bf 173} (2000) 103-153.}

 \end{thebibliography}
\end{document}